\documentclass[11pt,a4paper]{amsart}
\usepackage[T1]{fontenc}
\usepackage{lmodern,amsmath,amssymb,mathtools}
\usepackage[a4paper,margin=1in]{geometry}
\usepackage[hidelinks]{hyperref}
\numberwithin{equation}{section}
\allowdisplaybreaks[2]
\newtheorem{theorem}{Theorem}[section]
\newtheorem{lemma}[theorem]{Lemma}

\newtheorem{corollary}[theorem]{Corollary}
\theoremstyle{remark}
\newtheorem{remark}[theorem]{Remark}
\newcommand{\R}{\mathbb R}
\newcommand{\Hone}{H^1_0(\Omega)}
\newcommand{\dx}{\,\mathrm{d}x}
\newcommand{\dy}{\,\mathrm{d}y}
\newcommand{\dmu}[1]{|x|^{-#1}\dx}
\newcommand{\HS}[1]{2^*(#1)}
\newcommand{\norm}[1]{\lVert #1\rVert}
\newcommand{\diver}{\operatorname{div}}
\newcommand{\supp}{\operatorname{supp}}
\title[Sharp subcritical criterion and critical limits]{A Sharp Subcritical Existence Criterion and Critical Limits for Hardy--Sobolev Equations with Competing Nonlinearities}
\author{Zhi-Yun Tang}
\address{School of Mathematics and Statistics, Southwest University, Chongqing, 400715, People's Republic of China}
\email{tzy@swu.edu.cn}
\author{Ying Lv}
\address{School of Mathematics and Statistics, Southwest University, Chongqing, 400715, People's Republic of China}
\email{ly0904@swu.edu.cn}
\thanks{Corresponding author: Ying Lv.}
\hypersetup{pdfauthor={Zhi-Yun Tang, Ying Lv},pdftitle={A Sharp Subcritical Existence Criterion and Critical Limits for Hardy-Sobolev Equations with Competing Nonlinearities}}
\begin{document}
\begin{abstract}
We study positive finite-energy solutions of a Dirichlet problem with competing Hardy--Sobolev nonlinearities on a bounded connected $C^2$ domain in $\mathbb R^N$ ($N\geq3$) with a boundary singularity. The source and absorption terms are $|x|^{-s_2}u^{q-1}$ and $-\lambda|x|^{-s_1}u^{p-1}$, respectively, where $0\leq s_1<s_2<2$, $2<q<p<2^*(s_1)$ and $q\leq2^*(s_2)$, with $2^*(s)=2(N-s)/(N-2)$. For $q<2^*(s_2)$ we establish a sharp existence criterion: the problem has a positive solution for every $\lambda>0$ if and only if
\[
q>q_0(p):=2+\frac{2-s_2}{2-s_1}(p-2).
\]
No star-shapedness or curvature sign condition is required. The main tool is a localized Pohozaev identity based on a boundary-tangent vector field, which gives a bound on the Dirichlet energy of positive solutions in terms of their energy levels when $q>q_0(p)$. For fixed $2^*(s_2)<p<2^*(s_1)$, the estimate is uniform as $q$ approaches $2^*(s_2)$ from below and $\lambda$ varies in compact subsets of $(0,\infty)$. Under the additional conditions of negative mean curvature at the singular point and $p<2^*(s_1)-2/(N-2)$, subcritical approximation yields a positive critical solution for every $\lambda>0$, and the least energy among positive critical solutions is attained. For fixed $2^*(s_2)<p<2^*(s_1)$ and $\lambda>0$, if the critical problem has no positive solution, the maxima of positive subcritical solutions diverge and every choice of maximum point tends to the singular point as $q$ approaches $2^*(s_2)$. For solutions at the least mountain-pass level and least-energy positive solutions, we also determine the limiting energy and concentration measures.
\end{abstract}

\maketitle

\section{Introduction and main results}
Consider the following elliptic problem with a boundary singularity:
\begin{equation}\label{eq:problem}
\begin{cases}
-\Delta u=|x|^{-s_2}u^{q-1}-\lambda |x|^{-s_1}u^{p-1} &\text{in }\Omega,\\
u=0 &\text{on }\partial\Omega,
\end{cases}
\end{equation}
where $\Omega\subset\R^N$ is a bounded connected $C^2$ domain with $0\in\partial\Omega$, $N\geq3$, and $\lambda>0$. Throughout the paper, solutions of \eqref{eq:problem} are understood to be finite-energy weak solutions in $\Hone$. Let $\HS{s}=2(N-s)/(N-2)$ for $0\leq s<2$. Unless otherwise stated, we assume that
\begin{equation}\label{eq:assumptions}
0\leq s_1<s_2<2,\qquad 2<q<p<\HS{s_1},\qquad q<\HS{s_2}.
\end{equation}

The case $p>q=\HS{s_2}$ with $\lambda>0$ arises in the open problem of Li and Lin \cite[Remark 1.2]{Li2012}. In approaching this critical problem through subcritical equations, a basic question is to determine the exponents for which \eqref{eq:problem} has a positive solution for every $\lambda>0$. Set
\[
q_0(p)=2+\frac{2-s_2}{2-s_1}(p-2).
\]
Tang and Tang \cite[Theorem 1.1 and Lemma 2.2]{Tang2025a} obtained two positive solutions for small $\lambda$ and ruled out nonzero finite-energy weak solutions for large $\lambda$ when $2<q<q_0(p)$ on a bounded domain with $0\in\partial\Omega$. Their Theorem 1.3 establishes a global bifurcation threshold in a further restricted exponent range. Tang and Tang \cite[Theorem 1.1]{Tang2025b} proved that a positive solution exists for every $\lambda>0$ if $\Omega$ is star-shaped with respect to the origin and
\begin{equation}\label{eq:old}
\frac{\HS{s_1}-p}{p(p-2)}>\frac{\HS{s_2}-q}{q(q-2)}.
\end{equation}
Their proof is based on local minimization on the Nehari manifold. Under \eqref{eq:assumptions}, condition \eqref{eq:old} implies $q>q_0(p)$. For the equation with $-\Delta+1$ on $\R^N$, an existence result for every $\lambda>0$ under $q>q_0(p)$ was obtained in \cite[Theorem 1.3]{Tang2025c}.

Theorem \ref{thm:exist} gives a sharp answer to this question for bounded connected $C^2$ domains: under \eqref{eq:assumptions}, positive solutions exist for every $\lambda>0$ exactly when $q>q_0(p)$. The borderline $q=q_0(p)$ is excluded from the all-parameter existence range. Here sharpness concerns solvability for every positive parameter, not the determination of the set of admissible parameters for fixed exponents. In contrast to \cite[Theorem 1.1]{Tang2025b}, neither star-shapedness nor the additional exponent restriction \eqref{eq:old} is needed, although global $C^2$ boundary regularity is assumed.

The second question concerns the limit as the source exponent approaches its critical value. At that exponent, the geometry near the singular point plays a particular role. For a single critical Hardy--Sobolev term, Ghoussoub and Kang \cite[Theorem 1.1]{Ghoussoub2004} proved attainment of the best constant when $N\geq4$ and all the principal curvatures at the origin are negative. Ghoussoub and Robert \cite[Theorem 1.1]{Ghoussoub2006} replaced this assumption by negativity of the mean curvature and extended the result to $N=3$ in \cite[Theorem 1.1]{Ghoussoub2006a}. Lin and Wadade \cite[Theorems 1.1 and 1.2]{Lin2012a} studied related minimization problems with an additional unweighted Lebesgue-norm term.

For the one-critical case $q=\HS{s_2}$, Cerami, Zhong and Zou \cite[Theorems 1.5 and 1.6]{Cerami2015} obtained existence for small $\lambda$ and, under a further exponent restriction, for almost every $\lambda>0$ when the mean curvature at the origin is negative. Xu, Zhang and Zhong \cite[Theorem 1]{Xu2027} gave an explicit sufficient condition for small-parameter existence when the best Hardy--Sobolev constant of $\Omega$ is smaller than that of the half-space. He and Lv \cite[Theorem 1.1]{He2026} obtained existence for every $\lambda>0$ when the mean curvature at the origin is negative and $\HS{s_2}<p<\HS{s_1}-2/(N-2)$. For $\HS{s_1}-2/(N-2)\leq p<\HS{s_1}$, they proved that the set of parameters admitting positive solutions is bounded above, without a curvature sign assumption. In the $C^2$ setting their boundary regularity assumption is local near the origin, rather than global as in this paper. Other related work concerns different exponent regimes or solution classes. Tang and Tang \cite[Theorems 1.1 and 1.4]{Tang2025} obtained large-parameter nonexistence in the doubly critical case and a continuation result in the $L^1$ solution framework. Yang \cite[Theorem 1]{Yang2026} extended the latter continuation result to a wider range of exponents by constructing Poisson lifts of singular boundary data.

The exponent condition $q>q_0(p)$ enters the existence proof in two ways. It allows the construction of negative-energy endpoints by spatial scaling near the singular point. It also ensures positivity of the coefficient in the energy estimate obtained by combining a localized Pohozaev identity, based on a vector field tangent to the boundary, with the identity obtained by testing the equation with the solution itself. This estimate gives compactness of energy-bounded families at fixed subcritical exponents. Together with the monotonicity trick, it yields mountain-pass solutions for every $\lambda>0$. The necessity of the exponent condition follows from a pointwise estimate and the Hardy inequality.

The same energy estimate connects the subcritical existence criterion with the critical limits. When $\HS{s_2}<p<\HS{s_1}$, its constants can be chosen uniformly for $q$ near $\HS{s_2}$ and $\lambda$ in compact subsets of $(0,\infty)$. With the test-path construction of Cerami, Zhong and Zou \cite[Lemma 2.4(iv)]{Cerami2015} and the half-space Hardy--Sobolev inequality, this uniformity allows us to pass to the critical exponent at fixed $\lambda$ under the curvature and exponent conditions of Theorem \ref{thm:endpoint}. We also obtain attainment of the least energy among positive critical solutions and compactness of minimizing sequences. In the absence of a positive critical solution, we locate the blow-up of arbitrary positive subcritical solutions and determine the limiting energy and concentration measures for solutions at the least mountain-pass level and least-energy positive solutions.

We use the norm $\norm{u}=(\int_\Omega |\nabla u|^2\dx)^{1/2}$ on $\Hone$ and write $u^+=\max\{u,0\}$. The energy functional is
\begin{equation}\label{eq:energy}
J_\lambda(u)=\frac12\int_\Omega |\nabla u|^2\dx
 +\frac{\lambda}{p}\int_\Omega (u^+)^p\dmu{s_1}
 -\frac1q\int_\Omega (u^+)^q\dmu{s_2},\qquad u\in\Hone.
\end{equation}
For $e\in\Hone$, set
\[
\Gamma(e)=\{\gamma\in C([0,1],\Hone):\gamma(0)=0,\ \gamma(1)=e\}.
\]
Let $\Gamma_\lambda$ be the union of $\Gamma(e)$ over all nonnegative $e\in\Hone$ with $J_\lambda(e)<0$. When $\Gamma_\lambda\ne\varnothing$, we define the least mountain-pass level by
\begin{equation}\label{eq:globalmp}
c_\lambda=\inf_{\gamma\in\Gamma_\lambda}\max_{\tau\in[0,1]} J_\lambda(\gamma(\tau)).
\end{equation}

\begin{theorem}\label{thm:exist}
Suppose that \eqref{eq:assumptions} holds. Then problem \eqref{eq:problem} has a positive solution for every $\lambda>0$ if and only if $q>q_0(p)$.

If $q>q_0(p)$, then for every $\lambda>0$ the path class $\Gamma_\lambda$ is nonempty and there exists a positive solution $u_\lambda$ with $J_\lambda(u_\lambda)=c_\lambda>0$. Moreover, for every nonnegative $e\in\Hone$ with $J_\lambda(e)<0$, the fixed-endpoint level
\[
\inf_{\gamma\in\Gamma(e)}\max_{\tau\in[0,1]}J_\lambda(\gamma(\tau))
\]
is strictly positive and is attained by a positive solution. For each compact interval $\Lambda\subset(0,\infty)$, there are $m_\Lambda,M_\Lambda>0$ such that
\[
0<m_\Lambda\leq c_\lambda\leq M_\Lambda<\infty\quad(\lambda\in\Lambda).
\]
\end{theorem}

\begin{remark}\label{rem:range}
The additional exponent condition \eqref{eq:old} in \cite[Theorem 1.1]{Tang2025b} is strictly more restrictive than $q>q_0(p)$ under \eqref{eq:assumptions}. For fixed $p$, the right-hand side of \eqref{eq:old} is decreasing in $q\in(2,\HS{s_2})$. At $q=q_0(p)$ it exceeds the left-hand side by
\[
\frac{(s_2-s_1)(\HS{s_1}-p)}{p(2-s_1)q_0(p)}>0.
\]
This proves the implication. The converse fails: for $N=3$, $s_1=0$, $s_2=1$, $p=9/2$ and $q=33/10$, we have $q>q_0(p)$, whereas the difference between the left- and right-hand sides of \eqref{eq:old} is $-64/2145<0$.
\end{remark}

The following energy estimate is the main ingredient in the proof of sufficiency and in the passage to the critical exponent. We state it for both the subcritical and the one-critical cases.
\begin{theorem}\label{thm:estimate}
Suppose that $0\leq s_1<s_2<2$, $2<q<p<\HS{s_1}$, $q\leq\HS{s_2}$ and $q>q_0(p)$. There exists $c_0>0$, depending only on $N,s_1,s_2,p,q$, such that for every compact interval $\Lambda\subset(0,\infty)$ there is $K_\Lambda>0$ with
\begin{equation}\label{eq:coercivity}
J_\lambda(u)\geq c_0\left(\int_\Omega |\nabla u|^2\dx+\lambda\int_\Omega u^p\dmu{s_1}\right)-K_\Lambda
\end{equation}
for every positive solution $u$ of \eqref{eq:problem} with $\lambda\in\Lambda$. Moreover, for fixed $N,s_1,s_2,p$ and $\Lambda$, the constants can be chosen independently of $q\in[q_1,q_2]$ if $2<q_1\leq q_2\leq\HS{s_2}$, $q_2<p$ and $q_1>q_0(p)$.
\end{theorem}

\begin{corollary}\label{cor:compact}
Suppose that \eqref{eq:assumptions} holds and $q>q_0(p)$. Then for every compact interval $\Lambda\subset(0,\infty)$ and $M\in\R$, the set
\[
\mathcal S_{\Lambda,M}=\{(\lambda,u):\lambda\in\Lambda,\ u\text{ is a positive solution of }\eqref{eq:problem},\ J_\lambda(u)\leq M\}
\]
is compact in $\Lambda\times(\Hone\cap C^1(\overline\Omega))$, with the sum norm on the intersection space.
\end{corollary}

\begin{corollary}\label{cor:ground}
Under \eqref{eq:assumptions} and $q>q_0(p)$, for every $\lambda>0$ problem \eqref{eq:problem} has a positive solution $w_\lambda$ such that
\[
J_\lambda(w_\lambda)=\min\{J_\lambda(u):u\text{ is a positive solution of }\eqref{eq:problem}\}.
\]
\end{corollary}
\begin{remark}\label{rem:uniform}
Theorem \ref{thm:estimate} bounds positive solutions in terms of their energy levels; it does not assert coercivity of $J_\lambda$ on the whole space. Although the existence criterion is strictly subcritical, the energy estimate applies up to the critical source exponent. Indeed, there is a strict gap
\[
\HS{s_2}-q_0(p)=\frac{2-s_2}{2-s_1}(\HS{s_1}-p)>0.
\]
Consequently, when $\HS{s_2}<p<\HS{s_1}$, the exponent intervals in the energy estimate can include $\HS{s_2}$ while their lower endpoints remain strictly above $q_0(p)$. This is what allows the constants in the estimate to remain uniform in the subcritical approximation, including at the critical exponent itself.
\end{remark}

\begin{remark}\label{rem:ground}
Corollary \ref{cor:ground} minimizes the energy over positive solutions rather than over the Nehari set. The least-energy solution obtained in this way is not asserted to be distinct from the mountain-pass solution of Theorem \ref{thm:exist}.
\end{remark}

We next fix $p$ and consider the limit as $q$ approaches its critical value. Choose coordinates in which the tangent half-space at $0$ is $\R_+^N=\{x\in\R^N:x_N>0\}$ and denote by
\[
\mu_s(\R_+^N)=\inf_{0\ne w\in C_c^\infty(\R_+^N)}
\frac{\displaystyle\int_{\R_+^N}|\nabla w|^2\dx}
{\displaystyle\left(\int_{\R_+^N} |w|^{\HS{s}}\dmu{s}\right)^{2/\HS{s}}}
\]
the best Hardy--Sobolev constant of the half-space. Let $q_*:=\HS{s_2}$ and put
\[
c_+=\left(\frac12-\frac1{q_*}\right)\mu_{s_2}(\R_+^N)^{q_* /(q_*-2)}.
\]
When $q$ varies, we write $J_{\lambda,q}$ for the functional in \eqref{eq:energy} and $c_\lambda(q)$ for the level in \eqref{eq:globalmp}.

For the critical existence result, choose coordinates near $0$ such that
\[
\Omega\cap B_\rho=\{(x',x_N):x_N>h(x')\}\cap B_\rho,\qquad
h(0)=0,\quad \nabla h(0)=0.
\]
With the convention of \cite{Cerami2015}, $H(0)=\frac{1}{2(N-1)}\Delta h(0)$ for this graph representation; thus the negative mean-curvature condition is equivalent to $\Delta h(0)<0$.

\begin{theorem}\label{thm:endpoint}
Suppose that $0\leq s_1<s_2<2$, $q_*<p<\HS{s_1}-2/(N-2)$ and $H(0)<0$. Then for every $\lambda>0$ the following assertions hold.
\begin{enumerate}
\item[(i)] Let $q_n\in(2,q_*)$ with $q_n\to q_*$. For all sufficiently large $n$ there are positive solutions at the least mountain-pass levels $c_\lambda(q_n)$. Any sequence of such solutions has a subsequence converging weakly in $\Hone$ to a positive solution of \eqref{eq:problem} with $q=q_*$. In particular, the one-critical problem admits a positive solution.
\item[(ii)] The least energy
\[
m_\lambda=\inf\{J_{\lambda,q_*}(u):u\text{ is a positive solution of }\eqref{eq:problem}\text{ with }q=q_*\}
\]
is attained. Every minimizing sequence in this positive solution set has a strongly convergent subsequence in $\Hone$.
\item[(iii)] For every compact interval $\Lambda\subset(0,\infty)$, there are $\eta_\Lambda,C_\Lambda>0$ such that
\[
m_\lambda\leq c_+-\eta_\Lambda,\qquad \norm{w_\lambda}\leq C_\Lambda\quad(\lambda\in\Lambda)
\]
for every least-energy positive solution $w_\lambda$ of the one-critical problem.
\end{enumerate}
\end{theorem}
\begin{remark}\label{rem:criticalrange}
The exponent interval in Theorem \ref{thm:endpoint} is nonempty if and only if $s_2-s_1>1$, since
\[
\HS{s_1}-\frac{2}{N-2}-\HS{s_2}=\frac{2(s_2-s_1-1)}{N-2}.
\]
\end{remark}

\begin{remark}\label{rem:criticalmethod}
In the globally $C^2$ setting, part (i) of Theorem \ref{thm:endpoint} recovers the critical existence conclusion of \cite[Theorem 1.1(i)]{He2026}. Their argument passes from almost every parameter to every parameter at the critical exponent. Here we fix $\lambda$ and pass to the limit in the exponent, using the uniform estimate of Theorem \ref{thm:estimate} and the strict energy bound from the test path of \cite[Lemma 2.4(iv)]{Cerami2015}. Parts (ii) and (iii) describe the least energy in the positive solution set of the critical problem.
\end{remark}

We next consider a different limiting question: the behavior of subcritical solutions when the critical problem has no positive solution.

\begin{theorem}\label{thm:blowup}
Suppose that $0\leq s_1<s_2<2$ and $q_*<p<\HS{s_1}$. Fix $\lambda>0$ and let $q_n\in(2,q_*)$ with $q_n\to q_*$. Assume that problem \eqref{eq:problem} with $q=q_*$ has no positive solution. Then the following assertions hold.
\begin{enumerate}
\item[(i)] For any positive solutions $u_n$ of \eqref{eq:problem} with $q=q_n$,
\[
\max_{x\in\overline\Omega}u_n(x)=u_n(x_n)\to+\infty,\qquad x_n\to0,
\]
where $x_n$ is any maximum point of $u_n$.
\item[(ii)] For any such sequence of positive solutions,
\begin{equation}\label{eq:lowerlevel}
\liminf_{n\to\infty}J_{\lambda,q_n}(u_n)\geq c_+.
\end{equation}
If each $u_n$ is a solution at the least mountain-pass level or a least-energy positive solution, then
\begin{align}
\label{eq:criticallevel}
J_{\lambda,q_n}(u_n)&\to c_+,\\
\label{eq:criticalmass}
\int_\Omega|\nabla u_n|^2\dx&\to\mu_{s_2}(\R_+^N)^{q_* /(q_*-2)},\qquad
\int_\Omega u_n^p\dmu{s_1}\to0.
\end{align}
Under the same restriction on the sequence, $u_n\rightharpoonup0$ in $\Hone$ and $u_n\to0$ in $C^1$ on compact subsets of $\overline\Omega\setminus\{0\}$. In the sense of finite measures on $\overline\Omega$,
\begin{equation}\label{eq:measure}
|\nabla u_n|^2\dx,\quad |x|^{-s_2}u_n^{q_n}\dx
\ \rightharpoonup\ \mu_{s_2}(\R_+^N)^{q_* /(q_*-2)}\delta_0,
\end{equation}
where $\delta_0$ denotes the unit point mass at the origin.
\end{enumerate}
Solutions at the least mountain-pass level and least-energy positive solutions exist for all sufficiently large $n$ by Theorem \ref{thm:exist} and Corollary \ref{cor:ground}.
\end{theorem}
\begin{remark}\label{rem:blowupcases}
There are two useful settings in which the nonexistence hypothesis of Theorem \ref{thm:blowup} follows from earlier results. If $\Omega$ is star-shaped with respect to the origin, it holds for every $\lambda>0$ by \cite[Theorem 1.5]{Tang2025b}. In the exponent range $p\geq\HS{s_1}-2/(N-2)$, it holds for sufficiently large $\lambda$ by \cite[Theorem 1.1(ii)]{He2026}, whose local boundary hypotheses are satisfied in the present $C^2$ setting.
\end{remark}

\begin{remark}\label{rem:blowupcomparison}
For the star-shaped case, \cite[Theorem 1.6]{Tang2025b} obtains a subsequential boundary blow-up point for the positive solutions constructed under \eqref{eq:old}. Part (i) of Theorem \ref{thm:blowup} identifies the singular point as the limit of every choice of maximum point for arbitrary positive solution sequences under its hypotheses. The energy and measure limits in part (ii) for the two selected variational families are obtained from the half-space threshold $c_+$.
\end{remark}

\section{Preliminaries}\label{sec:prelim}
In what follows, $B_r$ denotes the ball of radius $r$ centered at the origin, and $C$ denotes a positive constant which may change from line to line.

\begin{lemma}\label{lem:embedding}
Let $0\leq s<2$ and $2\leq l<\HS{s}$. Then the embedding $\Hone\hookrightarrow L^l(\Omega,|x|^{-s}\dx)$ is compact.
\end{lemma}
\begin{proof}
By the Hardy and Sobolev inequalities for the zero extension of $u$, and H\"older's inequality, we have
\begin{align*}
\int_\Omega |u|^{\HS{s}}\dmu{s}
&\leq \left(\int_\Omega |x|^{-2}|u|^2\dx\right)^{s/2}
\left(\int_\Omega |u|^{\HS{0}}\dx\right)^{1-s/2}\\
&\leq C\norm{u}^{\HS{s}}.
\end{align*}
Here the case $s=0$ is just the Sobolev inequality. For $0<\varepsilon<1$, it follows from H\"older's inequality that
\begin{align*}
\int_{\Omega\cap B_\varepsilon}|u|^l\dmu{s}
&\leq\left(\int_\Omega|u|^{\HS{s}}\dmu{s}\right)^{l/\HS{s}}
\left(\int_{B_\varepsilon}|x|^{-s}\dx\right)^{1-l/\HS{s}}\\
&\leq C\varepsilon^{(N-s)(1-l/\HS{s})}\norm{u}^l.
\end{align*}
Suppose that $u_n\rightharpoonup u$ in $\Hone$. On $\Omega\setminus B_\varepsilon$ the weight is bounded, and the Rellich theorem gives convergence of $u_n$ to $u$ in the weighted $L^l$ norm on this set. Applying the last estimate to $u_n-u$ and then letting $\varepsilon\downarrow0$, we obtain the desired compactness.
\end{proof}

For $l>2$, the inequality
\[
\big|(a^+)^{l-1}-(b^+)^{l-1}\big|\leq(l-1)(a^++b^+)^{l-2}|a-b|\quad(a,b\in\R)
\]
and H\"older's inequality with respect to $|x|^{-s}\dx$ show that the map $v\mapsto (v^+)^{l-1}$ is continuous from $L^l(\Omega,|x|^{-s}\dx)$ into $L^{l/(l-1)}(\Omega,|x|^{-s}\dx)$. In particular, $J_\lambda\in C^1(\Hone,\R)$ also when $q=\HS{s_2}$, by the continuous Hardy--Sobolev embedding. If $J_\lambda'(u)=0$, testing with $u^-=\max\{-u,0\}$ gives $\norm{u^-}=0$; hence $u\geq0$ and $u$ satisfies the weak equation \eqref{eq:problem}. Under the strictly subcritical assumptions \eqref{eq:assumptions}, Lemma \ref{lem:embedding} implies that the nonlinear terms in its derivative map weakly convergent sequences in $\Hone$ to strongly convergent sequences in $H^{-1}(\Omega)$.

\begin{lemma}\label{lem:regularity}
Let $a,b>0$, $0\leq s_1<s_2<2$, $2<q<p<\HS{s_1}$ and $q\leq\HS{s_2}$. Suppose that $v\in\Hone$ is a nonzero nonnegative weak solution of
\[
-\Delta v=a|x|^{-s_2}v^{q-1}-b|x|^{-s_1}v^{p-1}.
\]
Then $v>0$ in $\Omega$, $v\in C^\infty(\Omega)\cap C^{1,\alpha}(\overline\Omega)$ for every $0<\alpha<\min\{1,q-s_2\}$, and $v$ belongs locally to $W^{2,m}$ and $C^{1,\alpha}$ on $\overline\Omega\setminus\{0\}$ for every $1<m<\infty$ and $0<\alpha<1$. If $q<\HS{s_2}$, then for each $0<\alpha<\min\{1,q-s_2\}$, the $C^{1,\alpha}(\overline\Omega)$ norm is bounded uniformly on families for which $a,b$ and $\norm{v}$ are bounded above (with fixed exponents and domain).
\end{lemma}
\begin{proof}
We first establish boundedness by a truncated Moser iteration. By the equation,
\[
-\Delta v\leq Vv,\qquad V=a|x|^{-s_2}v^{q-2}.
\]
For $k>0$ and $l\geq1$, put $v_k=\min\{v,k\}$ and $w=v v_k^{l-1}$, with $w=v$ when $l=1$. Testing the subsolution inequality with $v v_k^{2l-2}\in\Hone$ and using the chain rule gives
\[
\int_\Omega|\nabla w|^2\dx\leq l^2\int_\Omega \nabla v\cdot\nabla(v v_k^{2l-2})\dx
\leq l^2\int_\Omega Vw^2\dx.
\]
Indeed, on $\{v<k\}$ the two gradient densities before the factor $l^2$ on the right are $l^2v^{2l-2}|\nabla v|^2$ and $(2l-1)v^{2l-2}|\nabla v|^2$, respectively; on $\{v\geq k\}$ they coincide. Here the case $l=1$ uses the test function $v$.

For $q<\HS{s_2}$, the Sobolev embedding gives $v\in L^{\HS{0}}(\Omega)$, and $s_2/N+(q-2)/\HS{0}<2/N$. At $q=\HS{s_2}$, the strict inequality is lost, so we first improve the integrability of each fixed solution. Take $l=N/(N-2)$ and let $M>0$. Applying the preceding estimate, H\"older's inequality with respect to $|x|^{-s_2}\dx$ and the Hardy--Sobolev inequality,
\begin{align*}
\norm{w}^2&\leq al^2\int_\Omega v^{q-2}w^2\dmu{s_2}\\
&\leq Cal^2\left(\int_{\{v>M\}}v^q\dmu{s_2}\right)^{(q-2)/q}\norm{w}^2
+al^2M^{q-2}\int_\Omega w^2\dmu{s_2}.
\end{align*}
Since $v\in L^q(\Omega,|x|^{-s_2}\dx)$, choose $M$ so large that the coefficient of $\norm{w}^2$ on the right is less than $1/4$. As $|x|^{-s_2}\in L^m(\Omega)$ for some $m>N/2$, H\"older's inequality, interpolation and Young's inequality give
\[
al^2M^{q-2}\int_\Omega w^2\dmu{s_2}\leq \frac14\norm{w}^2+C\norm{w}_{L^2(\Omega)}^2.
\]
Thus $\norm{w}^2\leq C\norm{v}_{L^{2l}(\Omega)}^{2l}$ with $C$ independent of $k$, since $2l=\HS{0}$. The Sobolev inequality and monotone convergence as $k\to\infty$ give $v\in L^{l\HS{0}}(\Omega)$. Consequently, in both cases there is an exponent $r\geq\HS{0}$ with $v\in L^r(\Omega)$ and $s_2/N+(q-2)/r<2/N$.

With this integrability improvement, the same iteration applies in both cases. Choose $\ell_0>1$ such that $s_2/N<1/\ell_0<2/N-(q-2)/r$, and define $m_0$ by $1/m_0=1/\ell_0+(q-2)/r$. Then $m_0>N/2$ and
\[
\norm{V}_{L^{m_0}(\Omega)}\leq a\norm{|x|^{-s_2}}_{L^{\ell_0}(\Omega)}\norm{v}_{L^r(\Omega)}^{q-2}.
\]
For $q<\HS{s_2}$ one can take $r=\HS{0}$, so this bound is uniform if $a$ and $\norm{v}$ are bounded. We now use the preceding truncated test with arbitrary $l\geq1$. By the Sobolev and H\"older inequalities and interpolation, with $\theta=N/(2m_0)<1$, we obtain
\begin{align*}
\norm{w}_{L^{\HS{0}}}^2
&\leq Cl^2\norm{V}_{L^{m_0}}\norm{w}_{L^{2m_0/(m_0-1)}}^2\\
&\leq Cl^2\norm{V}_{L^{m_0}}\norm{w}_{L^2}^{2(1-\theta)}\norm{w}_{L^{\HS{0}}}^{2\theta}.
\end{align*}
Here and below the unweighted Lebesgue norms are taken on $\Omega$. Dividing when $w\ne0$ and then letting $k\to\infty$, for $v\in L^{2l}(\Omega)$ we have
\[
\norm{v}_{L^{l\HS{0}}}\leq (A l^2)^{1/(2l(1-\theta))}\norm{v}_{L^{2l}},\qquad A=C(1+\norm{V}_{L^{m_0}})\geq1.
\]
Taking $l_j=(N/(N-2))^j$ for $j=0,1,\ldots$ and using the convergence of $\sum_j l_j^{-1}$ and $\sum_j l_j^{-1}\log l_j$, we obtain by Moser iteration $\norm{v}_{L^\infty(\Omega)}\leq C\norm{v}_{L^2(\Omega)}$, where $C$ is bounded when $\norm{V}_{L^{m_0}}$ is bounded.

To obtain boundary regularity, we use the boundedness of $v$. The right-hand side of its equation belongs to $L^m(\Omega)$ for every $1<m<N/s_2$, with uniform bounds for fixed strictly subcritical exponents if $a,b$ and $\norm{v}$ are bounded. For $0<\gamma<\min\{1,2-s_2\}$, choose $N/(2-\gamma)<m<\min\{N,N/s_2\}$. By the Dirichlet solvability and estimate in \cite[Theorem 9.15 and Lemma 9.17]{Gilbarg2001}, the right-hand side determines a zero-trace $W^{2,m}$ solution. Since $m>N/2$, this solution belongs to $\Hone$ by Sobolev embedding, and hence equals $v$ by uniqueness of the weak Dirichlet problem. The Sobolev embedding then gives $v\in C^{0,\gamma}(\overline\Omega)$. Its continuous representative vanishes on $\partial\Omega$ by the trace condition, so that $v(x)\leq C|x|^\gamma$. On patches separated from $0$, the right-hand side is bounded, and interior and boundary regularity (see \cite[Theorem 9.13 and Lemma 9.16]{Gilbarg2001}) give the stated $W^{2,m}$ and $C^{1,\alpha}$ regularity away from $0$.

We improve the estimate at $0$ by using the vanishing of $v$ on the boundary, as in the regularity argument for a single Hardy--Sobolev term in \cite[Proposition A.1, Steps A.4 and A.5]{Ghoussoub2006a}. Suppose that $v(x)\leq C|x|^\gamma$ with $0<\gamma<1$ and put $t=s_2-(q-1)\gamma$. For $0<|x|<1$, since $s_1<s_2$ and $p>q$, we have
\[
|\Delta v|\leq C\big(|x|^{(q-1)\gamma-s_2}+|x|^{(p-1)\gamma-s_1}\big)\leq C|x|^{-t}.
\]
If $t\geq1$, then $1\leq t<2$ and
\[
2-t-\gamma=2-s_2+(q-2)\gamma\geq 2-s_2>0.
\]
By choosing $m<N/t$ sufficiently close to $N/t$, the global $W^{2,m}$ estimate and Sobolev embedding give $v\in C^{0,\gamma'}(\overline\Omega)$ with $\gamma'=\gamma+(2-s_2)/2<2-t\leq1$. Since $v(0)=0$, the bound $v(x)\leq C|x|^{\gamma'}$ follows. Repeating this argument finitely many times, we obtain an exponent for which $t<1$. We can then choose $m>N$ with $t m<N$ if $t>0$ (or any finite $m>N$ if $t\leq0$), and conclude that $v\in C^{1,\eta}(\overline\Omega)$ for some $\eta>0$.

Since $v(0)=0$ and $v\in C^1(\overline\Omega)$, a local boundary chart gives $v(x)\leq C|x|$ near $0$. Consequently, $|\Delta v|\leq C|x|^{q-1-s_2}$ there. Given $0<\alpha<\min\{1,q-s_2\}$, choose $m>N/(1-\alpha)$ such that $m(s_2+1-q)<N$ if $s_2+1-q>0$. This is possible by the restriction on $\alpha$. Applying the global $W^{2,m}$ estimate once more gives $v\in C^{1,\alpha}(\overline\Omega)$. For strictly subcritical $q$, all the constants in these estimates are uniform when $a,b$ and $\norm{v}$ are bounded above, with the exponents fixed. At $q=\HS{s_2}$ the choice of $M$ may depend on the solution.

Finally, to prove positivity and interior smoothness, we write
\[
-\Delta v+b|x|^{-s_1}v^{p-2}v=a|x|^{-s_2}v^{q-1}\geq0.
\]
The coefficient of $v$ on the left is nonnegative and locally bounded. Thus the strong maximum principle (see \cite[Theorem 8.19]{Gilbarg2001}) and connectedness of $\Omega$ give $v>0$ in $\Omega$. On each compact interior subset $v$ is bounded away from zero, so the nonlinearities are smooth functions of $(x,v)$. Interior elliptic regularity gives $v\in C^\infty(\Omega)$.
\end{proof}

\section{Uniform energy estimates for positive solutions}
In this section we allow $q\leq\HS{s_2}$ in \eqref{eq:assumptions}. We first establish bounds away from the singular point which do not depend on the energy of the solution. These bounds allow us to localize the Pohozaev identity near the origin without assuming star-shapedness. The resulting energy estimate also gives compactness for fixed strictly subcritical exponents.
\begin{lemma}\label{lem:away}
Let $\Lambda=[\lambda_-,\lambda_+]\subset(0,\infty)$. For each $r>0$ there exists $C_{r,\Lambda}>0$ such that every positive solution of \eqref{eq:problem} with $\lambda\in\Lambda$ satisfies
\begin{equation}\label{eq:away}
\sup_{\Omega\setminus B_r}u\leq C_{r,\Lambda},\qquad
\norm{u}_{C^{1,\alpha}(\overline\Omega\setminus B_{2r})}\leq C_{r,\Lambda,\alpha}
\end{equation}
for every fixed $\alpha\in(0,1)$, whenever the indicated sets are nonempty. The boundary norm is understood in local charts.
\end{lemma}
\begin{proof}
Set $R_* =1+\sup_{x\in\Omega}|x|$ and $a_0=\lambda_-R_*^{-s_1}/2$. Since $p>q$, there exists $M_r\geq0$ such that
\[
(r/2)^{-s_2}z^{q-1}-a_0z^{p-1}\leq M_r\quad(z\geq0).
\]
Hence for $x\in\Omega\setminus B_{r/2}$ and $\lambda\in\Lambda$,
\[
|x|^{-s_2}z^{q-1}-\lambda|x|^{-s_1}z^{p-1}\leq M_r-a_0z^{p-1},
\]
so that $-\Delta u+a_0u^{p-1}\leq M_r$ in $\Omega\setminus B_{r/2}$.

Fix $x_0\in\Omega\setminus B_r$ and let $R=r/4$, $\nu=2/(p-2)$. For $x\in B_R(x_0)$, put $d(x)=R^2-|x-x_0|^2$ and
\[
W(x)=K d(x)^{-\nu}.
\]
We have $\nu(p-1)=\nu+2$ and
\begin{align*}
\Delta(d^{-\nu})&=2\nu N d^{-\nu-1}+4\nu(\nu+1)|x-x_0|^2d^{-\nu-2}\\
&\leq [2\nu N+4\nu(\nu+1)]R^2d^{-\nu-2}=:C_Rd^{-\nu-2}.
\end{align*}
Consequently,
\[
-\Delta W+a_0W^{p-1}\geq(a_0K^{p-1}-C_RK)d^{-\nu-2}\geq M_r
\]
if $K$ is chosen so that $a_0K^{p-1}-C_RK\geq M_rR^{2(\nu+2)}$. This choice is independent of $u$ and $x_0$.

By Lemma \ref{lem:regularity}, the fixed solution $u$ is bounded on $\overline\Omega\cap\overline B_R(x_0)$, which is separated from the origin. Since $W$ tends uniformly to infinity as $|x-x_0|\uparrow R$, the function $(u-W)^+$, extended by zero outside the ball, vanishes in a neighborhood of its sphere and has zero trace on $\partial\Omega$. It is therefore an admissible nonnegative test function in $\Hone$. Subtracting the differential inequalities and testing with this function gives
\[
\int_{\Omega\cap B_R(x_0)}\!|\nabla(u-W)^+|^2\dx
+a_0\int_{\Omega\cap B_R(x_0)}\!(u^{p-1}-W^{p-1})(u-W)^+\dx\leq0.
\]
Both terms are nonnegative. Thus $u\leq W$ in the cap and, in particular, $u(x_0)\leq KR^{-2\nu}$. This proves the first estimate of \eqref{eq:away}. The right-hand side of \eqref{eq:problem} is now uniformly bounded on $\Omega\setminus B_r$. Applying the local and boundary $W^{2,m}$ estimates on smaller patches with $m>N/(1-\alpha)$ proves the second estimate in \eqref{eq:away}.
\end{proof}

To derive the energy estimate, we need to control the error in the Pohozaev identity associated with dilation. For nonnegative $u\in\Hone$ we set
\begin{equation}\label{eq:P}
P_\lambda(u)=\frac{N-2}{2}\int_\Omega|\nabla u|^2\dx
+\frac{\lambda(N-s_1)}p\int_\Omega u^p\dmu{s_1}
-\frac{N-s_2}q\int_\Omega u^q\dmu{s_2}.
\end{equation}
\begin{lemma}\label{lem:pohoz}
Let $\Lambda\subset(0,\infty)$ be a compact interval. For every $\varepsilon>0$ there exists $C_{\varepsilon,\Lambda}>0$ such that every positive solution of \eqref{eq:problem} with $\lambda\in\Lambda$ satisfies
\begin{equation}\label{eq:Psmall}
|P_\lambda(u)|\leq\varepsilon\left(\int_\Omega|\nabla u|^2\dx+\lambda\int_\Omega u^p\dmu{s_1}\right)+C_{\varepsilon,\Lambda}.
\end{equation}
\end{lemma}
\begin{proof}
To avoid a boundary term of indefinite sign, we construct a vector field tangent to the boundary near $0$. After a rotation, there are $r_0>0$ and a $C^2$ function $h$ with $h(0)=0$ and $\nabla h(0)=0$ such that
\[
\Omega\cap B_{4r_0}=\{(x',x_N):x_N>h(x')\}\cap B_{4r_0}.
\]
For $x\in B_{4r_0}$, let
\[
X(x)=(x',x_N-h(x')+x'\cdot\nabla h(x')).
\]
On $x_N=h(x')$, we have $X\cdot(-\nabla h,1)=0$, so $X\cdot\nu=0$ on the physical boundary, where $\nu$ is the outer unit normal. Moreover, $\partial_N X_N=1$ and $\partial_i X_N=\sum_{j=1}^{N-1}x_j\partial_{ij}h$ for $i<N$. It follows that
\begin{equation}\label{eq:X}
\begin{gathered}
\diver X=N,\qquad |X-x|\leq C|x|^2,\qquad |DX-I|\leq C|x|,\\
|X|\leq C|x|,\qquad \left|\frac{X\cdot x}{|x|^2}-1\right|\leq C|x|\quad (0<|x|<4r_0),
\end{gathered}
\end{equation}
where $I$ is the identity matrix and $C$ is independent of $u$, $\lambda$ and the radii below.

Fix $0<r<r_0$ and choose $\chi\in C_c^\infty(B_2)$ with $0\leq\chi\leq1$ and $\chi=1$ on $B_1$. Define $Y_r(x)=\chi(x/r)X(x)$ on the chart and extend it by zero to $\overline\Omega$. For $0<\eta<r/2$ put $Y_{r,\eta}(x)=(1-\chi(x/\eta))Y_r(x)$. Both fields are tangent to $\partial\Omega$ and $Y_{r,\eta}$ vanishes near $0$.

We use the following form of the identity in \cite[Proposition 1 and (4)]{Pucci1986}. Since the coefficients are singular at $0$, we first use fields which vanish near the origin. Let
\[
F_\lambda(x,z)=\frac{|x|^{-s_2}z^q}{q}-\frac{\lambda|x|^{-s_1}z^p}{p}\quad(z\geq0).
\]
For any $C^1$ vector field $Y$ tangent to $\partial\Omega$ and vanishing near $0$, the function $Y\cdot\nabla u$ belongs to $\Hone$. Indeed, Lemma \ref{lem:regularity} gives $W^{2,m}$ regularity on the support of $Y$, and $\nabla u=(\partial_\nu u)\nu$ on the boundary away from $0$. Testing \eqref{eq:problem} with $Y\cdot\nabla u$, we obtain
\begin{align*}
\int_\Omega \nabla u\cdot\nabla(Y\cdot\nabla u)\dx
&=\int_\Omega (DY\nabla u)\cdot\nabla u\dx
+\frac12\int_\Omega Y\cdot\nabla(|\nabla u|^2)\dx,\\
\int_\Omega \partial_z F_\lambda(x,u)Y\cdot\nabla u\dx
&=\int_\Omega Y\cdot\nabla(F_\lambda(x,u))\dx
-\int_\Omega Y\cdot\nabla_xF_\lambda(x,u)\dx.
\end{align*}
Here $\nabla_x F_\lambda$ denotes differentiation with $z$ fixed. Integrating by parts in the terms containing total gradients, the boundary integrals vanish since $Y\cdot\nu=0$. All terms are regular on the support of $Y$. Thus
\begin{equation}\label{eq:identity}
\begin{aligned}
0={}&\int_\Omega\left(\frac{\diver Y}{2}|\nabla u|^2-(DY\nabla u)\cdot\nabla u\right)\dx\\
&+\frac{\lambda}{p}\int_\Omega\left(\diver Y-s_1\frac{Y\cdot x}{|x|^2}\right)u^p\dmu{s_1}\\
&-\frac1q\int_\Omega\left(\diver Y-s_2\frac{Y\cdot x}{|x|^2}\right)u^q\dmu{s_2}.
\end{aligned}
\end{equation}
We next remove the cutoff at the singular point. Apply \eqref{eq:identity} with $Y=Y_{r,\eta}$. By \eqref{eq:X} and the cutoff construction, $|DY_{r,\eta}|+|Y_{r,\eta}|/|x|\leq C$ uniformly for $0<\eta<r/2$ and $0<r<r_0$. Moreover, $Y_{r,\eta}=Y_r$ outside $B_{2\eta}$. The difference between the integrals for these two fields is bounded by
\[
C\int_{\Omega\cap B_{2\eta}}\big(|\nabla u|^2+\lambda|x|^{-s_1}u^p+|x|^{-s_2}u^q\big)\dx,
\]
which tends to zero by the Hardy--Sobolev inequality and absolute continuity of the integrals. Thus \eqref{eq:identity} holds with $Y=Y_r$. No regularity of $u$ at $0$ is required.

It remains to compare the localized identity with \eqref{eq:P}. On $B_r$, $Y_r=X$. It follows from \eqref{eq:X} that the difference between the integrand in \eqref{eq:identity} and the integrand in \eqref{eq:P} is bounded in absolute value by
\[
Cr\big(|\nabla u|^2+\lambda|x|^{-s_1}u^p+|x|^{-s_2}u^q\big).
\]
On $\Omega\setminus B_r$ the integral of the difference is bounded by $C_{r,\Lambda}$, by Lemma \ref{lem:away} (applied with a smaller radius) and the bounds on the field. Consequently,
\[
|P_\lambda(u)|\leq Cr\left(\int_\Omega |\nabla u|^2\dx+\lambda\int_\Omega u^p\dmu{s_1}+\int_\Omega u^q\dmu{s_2}\right)+C_{r,\Lambda}.
\]
Testing the equation with $u$, we have the Nehari identity
\begin{equation}\label{eq:nehari}
\int_\Omega |\nabla u|^2\dx+\lambda\int_\Omega u^p\dmu{s_1}=\int_\Omega u^q\dmu{s_2}.
\end{equation}
Substitution into the preceding estimate gives
\[
|P_\lambda(u)|\leq 2Cr\left(\int_\Omega|\nabla u|^2\dx+\lambda\int_\Omega u^p\dmu{s_1}\right)+C_{r,\Lambda}.
\]
Since $C$ is independent of $r$, we may choose $r$ sufficiently small that $2Cr\leq\varepsilon$ and then fix it. This proves \eqref{eq:Psmall}.
\end{proof}

\begin{proof}[Proof of Theorem \ref{thm:estimate}]
We combine the Nehari identity with the Pohozaev estimate to obtain the coefficient involving $q-q_0(p)$. Using \eqref{eq:nehari} in \eqref{eq:energy} and \eqref{eq:P}, we have
\begin{align*}
J_\lambda(u)&=\frac{q-2}{2q}\int_\Omega u^q\dmu{s_2}-\frac{p-2}{2p}\lambda\int_\Omega u^p\dmu{s_1},\\
P_\lambda(u)&=\left(\frac{N-2}{2}-\frac{N-s_2}{q}\right)\int_\Omega u^q\dmu{s_2}\\
&\quad+\left(\frac{N-s_1}{p}-\frac{N-2}{2}\right)\lambda\int_\Omega u^p\dmu{s_1}.
\end{align*}
Consequently,
\begin{equation}\label{eq:combo}
\begin{aligned}
&\left(\frac{N-s_1}{p}-\frac{N-2}{2}\right)J_\lambda(u)+\frac{p-2}{2p}P_\lambda(u)\\
&\qquad=\frac{(2-s_1)(q-q_0(p))}{2pq}\left(\int_\Omega|\nabla u|^2\dx+\lambda\int_\Omega u^p\dmu{s_1}\right).
\end{aligned}
\end{equation}
Since $p<\HS{s_1}$ and $q>q_0(p)$, the coefficient of $J_\lambda$ and the coefficient on the right are both positive. Choose $\varepsilon>0$ sufficiently small that
\[
\frac{p-2}{2p}\varepsilon\leq\frac{(2-s_1)(q-q_0(p))}{4pq}.
\]
It follows from \eqref{eq:Psmall} and \eqref{eq:combo} that
\begin{align*}
\left(\frac{N-s_1}{p}-\frac{N-2}{2}\right)J_\lambda(u)
&\geq\frac{(2-s_1)(q-q_0(p))}{4pq}
\left(\int_\Omega|\nabla u|^2\dx+\lambda\int_\Omega u^p\dmu{s_1}\right)\\
&\quad-\frac{p-2}{2p}C_{\varepsilon,\Lambda}.
\end{align*}
Dividing by the positive coefficient of $J_\lambda(u)$ gives \eqref{eq:coercivity}.

If $q$ varies in $[q_1,q_2]$ as in the theorem, the constants in Lemma \ref{lem:away} can be chosen independently of $q$, since $z^{q-1}\leq1+z^{q_2-1}$ for $z\geq0$ and $p>q_2$. The geometric constants in \eqref{eq:X} do not depend on $q$, and the coefficients in \eqref{eq:identity} are bounded on $[q_1,q_2]$. The passage from $Y_{r,\eta}$ to $Y_r$ is made for each fixed solution. Thus the constant in \eqref{eq:Psmall} can also be chosen independently of $q$. Since $q_1>q_0(p)$, the coefficient on the right of \eqref{eq:combo} has a positive lower bound on $[q_1,q_2]$, so we may choose $\varepsilon$ and hence the constants in \eqref{eq:coercivity} independently of $q$. In particular, this does not require a uniform global regularity estimate at the critical exponent.
\end{proof}

\begin{proof}[Proof of Corollary \ref{cor:compact}]
Let $(\lambda_n,u_n)\in\mathcal S_{\Lambda,M}$. Theorem \ref{thm:estimate} gives boundedness of $(u_n)$ in $\Hone$. Passing to a subsequence, we may assume that $\lambda_n\to\lambda\in\Lambda$ and $u_n\rightharpoonup u$ in $\Hone$. By Lemma \ref{lem:embedding},
\[
u_n\to u\quad\text{in }L^p(\Omega,|x|^{-s_1}\dx)\text{ and }L^q(\Omega,|x|^{-s_2}\dx).
\]
In particular $u\geq0$, and by the nonlinear continuity following Lemma \ref{lem:embedding}, we may pass to the limit in the weak equation. Subtracting the equations for $u_n$ and $u$ and testing with $u_n-u$, we have
\begin{align*}
\norm{u_n-u}^2
&=\int_\Omega(u_n^{q-1}-u^{q-1})(u_n-u)\dmu{s_2}\\
&\quad-\int_\Omega(\lambda_nu_n^{p-1}-\lambda u^{p-1})(u_n-u)\dmu{s_1}\to0.
\end{align*}
Hence $u_n\to u$ in $\Hone$. From \eqref{eq:nehari} and the Hardy--Sobolev inequality,
\[
0<\norm{u_n}^2\leq\int_\Omega u_n^q\dmu{s_2}\leq C\norm{u_n}^q,
\]
whence $\norm{u_n}\geq C^{-1/(q-2)}$. The same lower bound holds for $\norm{u}$, so $u\ne0$. Lemma \ref{lem:regularity} gives $u>0$ in $\Omega$. The energy is continuous under this convergence and $\lambda_n\to\lambda$, so $J_\lambda(u)\leq M$. This gives compactness in $\Lambda\times\Hone$. Moreover, Lemma \ref{lem:regularity} with $a=1$ and $b=\lambda_n$ gives a uniform $C^{1,\alpha}(\overline\Omega)$ bound for a fixed $0<\alpha<\min\{1,q-s_2\}$. By the Arzel\`a--Ascoli theorem, after passing to a further subsequence, $u_n\to u$ in $C^1(\overline\Omega)$, where the limit is identified by the strong convergence in $\Hone$. This proves the assertion.
\end{proof}

\section{The sharp subcritical existence criterion}
We first prove the necessity of the strict exponent condition for solvability at every positive parameter. We then use spatial scaling and the preceding energy estimate to construct the mountain-pass solutions under $q>q_0(p)$.

\begin{proof}[Proof of necessity in Theorem \ref{thm:exist}]
Suppose that $q\leq q_0(p)$. Put
\[
\sigma=\frac{s_2(p-2)-s_1(q-2)}{p-q}.
\]
By $s_2>s_1$ and $p>q>2$, we have $\sigma>0$. Moreover,
\[
\sigma-2=\frac{(2-s_1)(q-q_0(p))}{p-q}\leq0.
\]
For $r>0$, the maximum of $r^{-s_2}z^{q-2}-\lambda r^{-s_1}z^{p-2}$ over $z\geq0$ is attained at the positive point satisfying
\[
z^{p-q}=\frac{q-2}{\lambda(p-2)}r^{s_1-s_2}.
\]
Substituting this value, we obtain
\[
r^{-s_2}z^{q-2}-\lambda r^{-s_1}z^{p-2}
\leq C\lambda^{-(q-2)/(p-q)}r^{-\sigma}\quad(z\geq0),
\]
where $C$ depends only on $p,q$. As $\Omega$ is bounded and $\sigma\leq2$, the Hardy inequality for the zero extension of any positive solution $u$ and \eqref{eq:nehari} give
\begin{align*}
0<\int_\Omega|\nabla u|^2\dx
&=\int_\Omega\left(|x|^{-s_2}u^{q-2}-\lambda|x|^{-s_1}u^{p-2}\right)u^2\dx\\
&\leq C\lambda^{-(q-2)/(p-q)}\int_\Omega |x|^{-2}u^2\dx\\
&\leq C\lambda^{-(q-2)/(p-q)}\int_\Omega|\nabla u|^2\dx,
\end{align*}
where $C$ is independent of $\lambda$ and $u$. This is impossible for sufficiently large $\lambda$, which proves the necessity.
\end{proof}

For the remainder of this section we assume \eqref{eq:assumptions} and $q>q_0(p)$. To apply Jeanjean's monotonicity trick \cite{Jeanjean1998,Jeanjean1999}, we put the parameter in front of the source term. For $t>0$, let
\[
I_t(v)=\frac12\int_\Omega|\nabla v|^2\dx+\frac1p\int_\Omega(v^+)^p\dmu{s_1}
-\frac tq\int_\Omega(v^+)^q\dmu{s_2},\qquad v\in\Hone.
\]
By the Hardy--Sobolev inequality and the continuity of the nonlinear terms established in Section \ref{sec:prelim}, $I_t\in C^1(\Hone,\R)$. Set
\begin{equation}\label{eq:rescale}
\lambda=t^{-(p-2)/(q-2)},\qquad u=t^{1/(q-2)}v.
\end{equation}
A direct calculation gives
\begin{equation}\label{eq:energyrescale}
J_\lambda(t^{1/(q-2)}v)=t^{2/(q-2)}I_t(v)\qquad(v\in\Hone).
\end{equation}
By differentiating \eqref{eq:energyrescale}, we see that critical points of $I_t$ are mapped to critical points of $J_\lambda$.

\begin{lemma}\label{lem:geometry}
For every interval $[a,b]\subset(0,\infty)$, there exist $\rho>0$ and a nonnegative $e\in\Hone$ with $\norm{e}>\rho$ such that
\[
I_t(v)\geq\frac14\norm{v}^2\quad\text{if }\norm{v}\leq\rho,\qquad I_t(e)<0\quad(t\in[a,b]).
\]
\end{lemma}
\begin{proof}
By the Hardy--Sobolev inequality,
\[
I_t(v)\geq\tfrac12\norm{v}^2-Cb\norm{v}^q\quad(t\in[a,b]).
\]
Since $q>2$, we can choose $\rho$ small enough that $Cb\rho^{q-2}\leq1/4$. The estimate holds on the closed ball and, in particular, gives the lower bound $\eta=\rho^2/4$ on its boundary.

To obtain a common negative endpoint, we use concentration at the origin. In the coordinates used in the proof of Lemma \ref{lem:pohoz}, take a nonzero nonnegative $w\in C_c^\infty(\R_+^N)$. Since $h(x')=O(|x'|^2)$, $r\supp w\subset\Omega$ for all sufficiently small $r>0$. Set
\begin{equation}\label{eq:spatial}
\beta=\frac{2-s_2}{q-2},\qquad \delta=2-s_1-\beta(p-2)=\frac{(2-s_1)(q-q_0(p))}{q-2}>0.
\end{equation}
For $L>0$, let $e_{r,L}(x)=Lr^{-\beta}w(x/r)$, extended by zero. Since $N-s_2-\beta q=N-2-2\beta$ and $N-s_1-\beta p=N-2-2\beta+\delta$, a change of variables gives
\begin{align*}
I_t(e_{r,L})=r^{N-2-2\beta}\bigg(&\frac{L^2}{2}\int_{\R_+^N}|\nabla w|^2\dy-\frac{tL^q}{q}\int_{\R_+^N}|y|^{-s_2}w^q\dy\\
&+\frac{r^\delta L^p}{p}\int_{\R_+^N}|y|^{-s_1}w^p\dy\bigg).
\end{align*}
Choose $L$ sufficiently large that the sum of the first two terms in parentheses is negative at $t=a$. Keeping $L$ fixed, we next choose $r$ sufficiently small that $I_a(e_{r,L})<0$, using $\delta>0$. Moreover, $q<\HS{s_2}$ implies $\beta>(N-2)/2$, so
\[
\norm{e_{r,L}}=Lr^{(N-2-2\beta)/2}\norm{\nabla w}_{L^2(\R_+^N)}\to\infty\quad(r\downarrow0).
\]
Taking $r$ smaller if necessary and putting $e=e_{r,L}$, we obtain $\norm e>\rho$ and $I_t(e)\leq I_a(e)<0$ for $t\in[a,b]$.
\end{proof}

The spatial scaling in the preceding proof also gives a relation between equations on dilated domains. Let $\beta,\delta$ be given by \eqref{eq:spatial} and, for $r>0$, set $\Omega_r=r^{-1}\Omega$. For a positive solution $u$ of \eqref{eq:problem}, $v(y)=r^\beta u(ry)$ satisfies
\[
\begin{cases}
-\Delta v=|y|^{-s_2}v^{q-1}-\lambda r^\delta|y|^{-s_1}v^{p-1}&\text{in }\Omega_r,\\
v=0&\text{on }\partial\Omega_r.
\end{cases}
\]
Conversely, a positive solution of the latter problem gives a solution of \eqref{eq:problem} by $u(x)=r^{-\beta}v(x/r)$. Writing $J^U_\lambda$ for the energy integrated over $U$, one has
\[
J^\Omega_\lambda(u)=r^{N-2-2\beta}J^{\Omega_r}_{\lambda r^\delta}(v).
\]
Indeed, the Dirichlet, source and absorbing integrals acquire the factors $r^{N-2-2\beta}$, $r^{N-s_2-\beta q}$ and $r^{N-s_1-\beta p}$, respectively. The identities following \eqref{eq:spatial} give the energy formula, and the same change of variables in the weak equation gives the transformed problem. Dilation preserves the zero-trace Sobolev spaces for each fixed $r>0$. The formulas remain valid without the sign restriction on $\delta$ used in Lemma \ref{lem:geometry}.

We now fix an interval $[a,b]\subset(0,\infty)$ and any nonnegative common negative-energy endpoint $e$. By the small-ball estimate, $\norm{e}>\rho$. Let
\[
c_e(t)=\inf_{\gamma\in\Gamma(e)}\max_{\tau\in[0,1]}I_t(\gamma(\tau)).
\]
Every path in $\Gamma(e)$ crosses the sphere of radius $\rho$. Using the straight path $\tau\mapsto\tau e$ for an upper bound, we have
\begin{equation}\label{eq:levels}
0<\eta:=\rho^2/4\leq c_e(t)\leq\frac12\int_\Omega|\nabla e|^2\dx+\frac1p\int_\Omega e^p\dmu{s_1}=:M_e\quad(t\in[a,b]).
\end{equation}
\begin{lemma}\label{lem:left}
The function $c_e$ is nonincreasing on $[a,b]$ and left-continuous on $(a,b]$.
\end{lemma}
\begin{proof}
The monotonicity follows from $I_{t_2}\leq I_{t_1}$ for $t_2\geq t_1$ and the common path class. As in \cite[Lemma 2.3]{Jeanjean1999}, take $t_n\uparrow t_0\in(a,b]$ and, for $\xi>0$, choose $\gamma\in\Gamma(e)$ with $\max_\tau I_{t_0}(\gamma(\tau))\leq c_e(t_0)+\xi$. Then
\[
c_e(t_0)\leq c_e(t_n)\leq c_e(t_0)+\xi+\frac{t_0-t_n}{q}\max_{\tau\in[0,1]}\int_\Omega(\gamma(\tau)^+)^q\dmu{s_2}.
\]
The maximum is finite by continuity and compactness of the path image. Letting $n\to\infty$ and then $\xi\downarrow0$ proves the assertion.
\end{proof}

\begin{lemma}\label{lem:ae}
For almost every $t\in[a,b]$, $I_t$ has a positive critical point $v_t$ with $I_t(v_t)=c_e(t)$.
\end{lemma}
\begin{proof}
The parameter-independent part of $I_t$ is bounded below by $\norm{v}^2/2$, and the integral multiplying $-t$ is nonnegative. By Lemma \ref{lem:geometry} and \eqref{eq:levels}, the assumptions of \cite[Theorem 1.1]{Jeanjean1999} hold for the common path class $\Gamma(e)$. Hence for almost every $t$ there exists a bounded sequence $(v_n)$ in $\Hone$ such that $I_t(v_n)\to c_e(t)$ and $I_t'(v_n)\to0$ in $H^{-1}(\Omega)$.

Passing to a subsequence, $v_n\rightharpoonup v$ in $\Hone$ and $v_n\to v$ in the two weighted spaces by Lemma \ref{lem:embedding}. The continuity of the nonlinear terms gives
\[
\norm{v_n-v}^2=\langle I_t'(v_n),v_n-v\rangle+o(1)\to0,
\]
where the gradient term involving $v$ tends to zero by weak convergence. Thus $I_t'(v)=0$ and $I_t(v)=c_e(t)\geq\eta$. Testing the equation with $v^-=\max\{-v,0\}$, we have $-\norm{v^-}^2=0$. Therefore $v\geq0$, $v\ne0$, and Lemma \ref{lem:regularity} gives $v>0$ in $\Omega$.
\end{proof}

\begin{proof}[Proof of sufficiency in Theorem \ref{thm:exist}]
We show that the level $c_e(t)$ is attained for every $t\in[a,b]$. Since $I_a(e)<0$, there exists $a'\in(0,a)$ with $I_{a'}(e)<0$. The same endpoint and the same small-sphere estimate are valid on $[a',b]$. We may therefore use Lemmas \ref{lem:left} and \ref{lem:ae} on this larger interval, with the same path class.

For any $t_0\in[a,b]$, choose $t_n\uparrow t_0$ from the full-measure set in Lemma \ref{lem:ae} in $(a',t_0)$ and positive critical points $v_n$ such that $I_{t_n}(v_n)=c_e(t_n)$. Put
\[
\lambda_n=t_n^{-(p-2)/(q-2)},\qquad u_n=t_n^{1/(q-2)}v_n.
\]
By \eqref{eq:rescale}--\eqref{eq:energyrescale}, the functions $u_n$ are positive solutions at $\lambda_n$, their energies are uniformly bounded above, and $\lambda_n\to\bar\lambda=t_0^{-(p-2)/(q-2)}>0$. Corollary \ref{cor:compact} gives a subsequence such that $u_n\to u_0$ strongly in $\Hone$, where $u_0$ is a positive solution at $\bar\lambda$. Let $v_0=t_0^{-1/(q-2)}u_0$. Then $v_n\to v_0$ in $\Hone$ and $I_{t_0}'(v_0)=0$. By Lemma \ref{lem:left},
\[
I_{t_0}(v_0)=\lim_{n\to\infty}I_{t_n}(v_n)=\lim_{n\to\infty}c_e(t_n)=c_e(t_0)>0.
\]
Thus, for this endpoint, there is a positive critical point $v_{t_0}$ at $c_e(t_0)$.

Given $\lambda>0$ and any nonnegative $e_\lambda$ with $J_\lambda(e_\lambda)<0$, take $t_0=\lambda^{-(q-2)/(p-2)}$ and $e=t_0^{-1/(q-2)}e_\lambda$. By \eqref{eq:energyrescale}, $I_{t_0}(e)<0$. By continuity, the same inequality holds on a compact interval $[a,b]$ containing $t_0$ in its interior. The small-ball estimate of Lemma \ref{lem:geometry} does not depend on the choice of $e$, so the preceding argument applies to this endpoint. Multiplication by $t_0^{1/(q-2)}$ maps $\Gamma(e)$ bijectively onto $\Gamma(e_\lambda)$. Thus $t_0^{1/(q-2)}v_{t_0}$ is a positive solution at the fixed-endpoint mountain-pass level of $J_\lambda$.

By Lemma \ref{lem:geometry} and \eqref{eq:energyrescale}, $\Gamma_\lambda\ne\varnothing$. For $\lambda$ in a compact interval $\Lambda$, choose an interval $[a,b]$ containing the corresponding $t$-values and a common endpoint $e$ from Lemma \ref{lem:geometry}. The small-ball estimate applies to every path in $\Gamma_\lambda$ after the inverse amplitude transformation, whereas the straight path to the common endpoint gives an upper bound. Hence
\[
0<\eta a^{2/(q-2)}\leq c_\lambda\leq M_eb^{2/(q-2)}<\infty.
\]

To show that the infimum over endpoints is attained as a critical value, fix $\lambda$ and choose a sequence of nonnegative endpoints $e_j$ with $J_\lambda(e_j)<0$ such that their fixed-endpoint levels converge to $c_\lambda$. By the first part of the proof, there are positive solutions $w_j$ at these levels. Their energies are bounded above. Corollary \ref{cor:compact} gives a subsequence converging in $\Hone$ to a positive solution $u_\lambda$. By continuity of $J_\lambda$, $J_\lambda(u_\lambda)=c_\lambda$.
\end{proof}

\begin{proof}[Proof of Corollary \ref{cor:ground}]
Fix $\lambda>0$. Theorem \ref{thm:exist} gives a positive solution, and Theorem \ref{thm:estimate} bounds the energies of all positive solutions below. Take a sequence of positive solutions $(w_n)$ whose energies tend to their infimum. Corollary \ref{cor:compact}, applied at the fixed parameter $\lambda$, gives a subsequence converging in $\Hone$ to a positive solution $w_\lambda$. By continuity of $J_\lambda$, $w_\lambda$ attains the infimum.
\end{proof}

\section{Critical limits}
In this section assume $q_*=\HS{s_2}<p<\HS{s_1}$. The gap in Remark \ref{rem:uniform} allows us to fix $q_1\in(q_0(p),q_*)$. Thus the subcritical existence criterion applies throughout $[q_1,q_*)$, and Theorem \ref{thm:estimate} gives a uniform energy estimate for $q\in[q_1,q_*]$ and $\lambda$ in compact subsets of $(0,\infty)$.

\subsection{Critical weak limits}
The next lemma gives the limiting equation and the energy bound for a zero weak limit. It will be used in both the critical existence and blow-up arguments.

\begin{lemma}\label{lem:criticalweak}
Let $q_n\in[q_1,q_*]$ with $q_n\to q_*$, $\lambda_n\to\lambda>0$, and let $u_n$ be positive solutions of \eqref{eq:problem} with exponents $q_n$ and parameters $\lambda_n$. Suppose that $u_n\rightharpoonup u$ in $\Hone$. Then $u$ is a nonnegative solution at $(\lambda,q_*)$, and
\begin{equation}\label{eq:criticallsc}
J_{\lambda,q_*}(u)\leq \liminf_{n\to\infty}J_{\lambda_n,q_n}(u_n).
\end{equation}
If $u=0$, then
\begin{equation}\label{eq:zerothreshold}
\liminf_{n\to\infty}J_{\lambda_n,q_n}(u_n)\geq c_+.
\end{equation}
\end{lemma}
\begin{proof}
By Lemma \ref{lem:embedding}, $u_n\to u$ in $L^p(\Omega,|x|^{-s_1}\dx)$ and in $L^2(\Omega)$. By Lemma \ref{lem:away} and the uniform estimates in the proof of Theorem \ref{thm:estimate}, the sequence is bounded in $C^{1,\alpha}$ locally on $\overline\Omega\setminus\{0\}$, for every $0<\alpha<1$. Thus $u_n\to u$ in $C^1$ on each compact subset of $\overline\Omega\setminus\{0\}$: every locally convergent subsequence has the same limit by the weak convergence. Passing to the limit in the weak equations against test functions in $C_c^\infty(\Omega)$ gives the equation at $(\lambda,q_*)$. Both nonlinear terms for $u$ define elements of $H^{-1}(\Omega)$ by the Hardy--Sobolev inequality. The weak equation therefore holds for every test function in $\Hone$ by density. In particular, $u\geq0$ and, if $u\ne0$, then $u>0$ in $\Omega$ by Lemma \ref{lem:regularity}.

Using \eqref{eq:nehari}, for both $u_n$ and the weak limit $u$, we obtain
\begin{align*}
J_{\lambda_n,q_n}(u_n)&=\left(\frac12-\frac1{q_n}\right)\norm{u_n}^2
-\lambda_n\left(\frac1{q_n}-\frac1p\right)\int_\Omega u_n^p\dmu{s_1},\\
J_{\lambda,q_*}(u)&=\left(\frac12-\frac1{q_*}\right)\norm{u}^2
-\lambda\left(\frac1{q_*}-\frac1p\right)\int_\Omega u^p\dmu{s_1}.
\end{align*}
The strong convergence in the weighted $L^p$ space and the weak lower semicontinuity of the Dirichlet integral give \eqref{eq:criticallsc}.

By the Hardy--Sobolev inequality and H\"older's inequality with respect to the finite measure $|x|^{-s_2}\dx$, there is $C\geq1$ independent of $n$ such that
\begin{equation}\label{eq:uniformnonzero}
0<\norm{u_n}^2\leq\int_\Omega u_n^{q_n}\dmu{s_2}\leq C\norm{u_n}^{q_n}.
\end{equation}
Hence $\norm{u_n}\geq C^{-1/(q_1-2)}>0$. Suppose now that $u=0$. Along a subsequence realizing the lower limit in \eqref{eq:zerothreshold}, we may assume that $\norm{u_n}^2\to\ell>0$. Since the absorbing integral tends to zero, the Nehari and energy identities give
\begin{equation}\label{eq:massbalance}
\int_\Omega u_n^{q_n}\dmu{s_2}\to\ell,\qquad
J_{\lambda_n,q_n}(u_n)\to\left(\frac12-\frac1{q_*}\right)\ell.
\end{equation}
We show that $\ell\geq\mu_{s_2}(\R_+^N)^{q_* /(q_*-2)}$. Fix a small $r>0$ and take $\chi_r\in C_c^\infty(B_{2r})$ with $0\leq\chi_r\leq1$ and $\chi_r=1$ on $B_r$. For $v\in\Hone$ supported in $B_{2r}$, the boundary chart in Lemma \ref{lem:pohoz} gives
\begin{equation}\label{eq:localHS}
\mu_{s_2}(\R_+^N)\left(\int_\Omega |v|^{q_*}\dmu{s_2}\right)^{2/q_*}
\leq(1+Cr)\int_\Omega |\nabla v|^2\dx,
\end{equation}
where $C$ is independent of small $r$. Indeed, the change of variables $x=(y',y_N+h(y'))$ has Jacobian one, its derivative is $I+O(r)$, and $|x|/|y|=1+O(r)$ on the support of $v$. The pullback of $v$ has zero trace on the flat boundary. Extending it by zero in $\R_+^N$ outside the chart gives a function in the completion of $C_c^\infty(\R_+^N)$ in the Dirichlet norm. The definition of $\mu_{s_2}(\R_+^N)$, the Hardy--Sobolev inequality and the change of variables give \eqref{eq:localHS}.

For each fixed $r$, since $u_n\to0$ in $L^2(\Omega)$, we have
\begin{align*}
\int_\Omega|\nabla(\chi_ru_n)|^2\dx
&=\int_\Omega\chi_r^2|\nabla u_n|^2\dx
+2\int_\Omega\chi_ru_n\nabla\chi_r\cdot\nabla u_n\dx
+\int_\Omega u_n^2|\nabla\chi_r|^2\dx\\
&\leq\norm{u_n}^2+o(1).
\end{align*}
The convergence to zero away from $0$ gives $\int_{\Omega\cap B_r}u_n^{q_n}\dmu{s_2}\to\ell$ by \eqref{eq:massbalance}. Furthermore,
\[
\int_{\Omega\cap B_r}u_n^{q_n}\dmu{s_2}
\leq\left(\int_\Omega|\chi_ru_n|^{q_*}\dmu{s_2}\right)^{q_n/q_*}
\left(\int_{\Omega\cap B_r}|x|^{-s_2}\dx\right)^{1-q_n/q_*}.
\]
Applying \eqref{eq:localHS} to $\chi_ru_n$, letting $n\to\infty$ and using $q_n\to q_*$, we obtain
\[
\mu_{s_2}(\R_+^N)\ell^{2/q_*}\leq(1+Cr)\ell.
\]
Letting $r\downarrow0$ and using $\ell>0$, we conclude that $\ell\geq\mu_{s_2}(\R_+^N)^{q_* /(q_*-2)}$. It follows from \eqref{eq:massbalance} that \eqref{eq:zerothreshold} holds.
\end{proof}

\subsection{Existence at the critical exponent}
The curvature and exponent conditions in Theorem \ref{thm:endpoint} give a test path below the half-space energy $c_+$. We use this strict inequality to exclude a zero weak limit in the subcritical approximation.

\begin{proof}[Proof of Theorem \ref{thm:endpoint}]
Fix $\Lambda=[\lambda_-,\lambda_+]\subset(0,\infty)$. Since $\Delta h(0)<0$ and
\[
p-1<\frac{N-2s_1}{N-2},
\]
by \cite[Lemma 2.4(iv)]{Cerami2015}, applied with their parameter $-\lambda_+$ and their exponent $p-1$, there exists a nonzero nonnegative $e\in\Hone$ such that
\[
J_{\lambda_+,q_*}(e)<0,\qquad\max_{0\leq\tau\leq1}J_{\lambda_+,q_*}(\tau e)<c_+.
\]
Since $J_{\lambda,q_*}\leq J_{\lambda_+,q_*}$ for $\lambda\in\Lambda$, the same endpoint works for all $\lambda\in\Lambda$. For this fixed $e$, the integrands $(\tau e)^q/q$, $q\in[q_1,q_*]$, $0\leq\tau\leq1$, are bounded by $C(1+e^{q_*})$ in the weighted space with measure $|x|^{-s_2}\dx$. Their dependence on $(\tau,q)$ is continuous, including at $\tau=0$. By dominated convergence and compactness of the parameter interval for $\tau$, the convergence of $J_{\lambda,q}(\tau e)$ to $J_{\lambda,q_*}(\tau e)$ as $q\to q_*^-$ is uniform in $\tau$ and $\lambda\in\Lambda$ (the other two terms are independent of $q$). Consequently there exist $\eta_\Lambda>0$ and $q_\Lambda\in(q_1,q_*)$ such that
\begin{equation}\label{eq:endpointgap}
J_{\lambda,q}(e)<0,\qquad\max_{0\leq\tau\leq1}J_{\lambda,q}(\tau e)\leq c_+-\eta_\Lambda
\quad(\lambda\in\Lambda,\ q_\Lambda\leq q\leq q_*).
\end{equation}

Fix $\lambda\in\Lambda$ and let $q_n<q_*$, $q_n\to q_*$. By Theorem \ref{thm:exist} and \eqref{eq:endpointgap}, for sufficiently large $n$ there are positive solutions $u_n$ at the least mountain-pass level, and every such solution satisfies
\[
0<J_{\lambda,q_n}(u_n)=c_\lambda(q_n)\leq\max_{0\leq\tau\leq1}J_{\lambda,q_n}(\tau e)\leq c_+-\eta_\Lambda.
\]
It follows from Theorem \ref{thm:estimate} that $(u_n)$ is bounded in $\Hone$. Take a subsequence such that $u_n\rightharpoonup u$ in $\Hone$. By \eqref{eq:zerothreshold}, $u$ cannot be zero. Lemma \ref{lem:criticalweak} now gives a positive solution at $q=q_*$ with
\begin{equation}\label{eq:endpointenergy}
J_{\lambda,q_*}(u)\leq\liminf_{n\to\infty}J_{\lambda,q_n}(u_n)\leq c_+-\eta_\Lambda.
\end{equation}
This proves (i).

By Theorem \ref{thm:estimate} and \eqref{eq:endpointenergy}, $-\infty<m_\lambda\leq c_+-\eta_\Lambda$. Let $(w_j)$ be a minimizing sequence in the positive solution set of the critical equation. Again by Theorem \ref{thm:estimate}, after passing to a subsequence, $w_j\rightharpoonup w$ in $\Hone$. Applying Lemma \ref{lem:criticalweak} with $q_n=q_*$ and constant parameter $\lambda$, we have $w\ne0$, because $m_\lambda<c_+$. Thus $w$ is a positive solution and
\[
m_\lambda\leq J_{\lambda,q_*}(w)\leq\liminf_{j\to\infty}J_{\lambda,q_*}(w_j)=m_\lambda.
\]
In particular, $w$ attains $m_\lambda$. Moreover, $w_j\to w$ in $L^p(\Omega,|x|^{-s_1}\dx)$ by Lemma \ref{lem:embedding}, and
\begin{align*}
\left(\frac12-\frac1{q_*}\right)\norm{w_j}^2
&=J_{\lambda,q_*}(w_j)+\lambda\left(\frac1{q_*}-\frac1p\right)\int_\Omega w_j^p\dmu{s_1}\\
&\to J_{\lambda,q_*}(w)+\lambda\left(\frac1{q_*}-\frac1p\right)\int_\Omega w^p\dmu{s_1}\\
&=\left(\frac12-\frac1{q_*}\right)\norm{w}^2.
\end{align*}
Hence $w_j\to w$ in $\Hone$. This proves (ii). Finally, applying the uniform estimate in Theorem \ref{thm:estimate} to any minimizer at $q=q_*$ and using \eqref{eq:endpointenergy} gives $\norm{w_\lambda}\leq C_\Lambda$ for $\lambda\in\Lambda$. This proves (iii).
\end{proof}

\subsection{Blow-up and energy concentration}
We now study the limit when the critical problem has no positive solution. For the equation with one critical Hardy--Sobolev source and a linear lower-order term, Ghoussoub and Robert \cite[Theorems 1.3 and 1.4]{Ghoussoub2006} studied minimal-energy blow-up and existence, and \cite[Theorem 1.3]{Ghoussoub2006a} established compactness of energy-bounded subcritical families under appropriate curvature assumptions. For the present problem, we use Lemma \ref{lem:criticalweak} to identify the limiting energy.

\begin{proof}[Proof of Theorem \ref{thm:blowup}(i)]
Discard finitely many terms so that $q_n\geq q_1$. Each $u_n$ is continuous on $\overline\Omega$ by Lemma \ref{lem:regularity}, and its maximum is attained in $\Omega$. Suppose, contrary to the assertion, that there exists a subsequence bounded in $L^\infty(\Omega)$. Then by \eqref{eq:nehari},
\[
\norm{u_n}^2\leq\int_\Omega u_n^{q_n}\dmu{s_2}\leq C.
\]
Take a further subsequence converging weakly in $\Hone$ and almost everywhere to a nonnegative function $u$. By Lemma \ref{lem:criticalweak} and the nonexistence assumption, $u=0$. Since $u_n$ is bounded in $L^\infty$, the dominated convergence theorem gives $\int_\Omega u_n^{q_n}\dmu{s_2}\to0$, in contradiction with the uniform nonzero bound \eqref{eq:uniformnonzero} (whose derivation uses only the Nehari and Hardy--Sobolev inequalities). Thus $\max_{\overline\Omega}u_n\to+\infty$. For every fixed $r>0$, Lemma \ref{lem:away}, with the uniform constants described in the proof of Theorem \ref{thm:estimate}, gives a bound on $u_n$ in $\Omega\setminus B_r$ independent of $n$. Therefore every maximum point $x_n$ belongs to $B_r$ for sufficiently large $n$. Since $r>0$ is arbitrary, $x_n\to0$.
\end{proof}

\begin{proof}[Proof of Theorem \ref{thm:blowup}(ii)]
We first give an upper bound for $c_\lambda(q)$ as $q\to q_*^-$, without any curvature assumption. Take a nonzero nonnegative $w\in C_c^\infty(\R_+^N)$ and set
\[
w_r(x)=r^{-(N-2)/2}w(x/r)
\]
for sufficiently small $r>0$, extending it by zero to $\Omega$. The function belongs to $\Hone$, since $\R_+^N$ is the tangent half-space at $0$ in the chosen coordinates. Put $\delta_*=N-s_1-(N-2)p/2>0$. At the critical exponent a change of variables gives
\begin{align*}
J_{\lambda,q_*}(\tau w_r)
={}&\frac{\tau^2}{2}\int_{\R_+^N}|\nabla w|^2\dy-\frac{\tau^{q_*}}{q_*}\int_{\R_+^N} |y|^{-s_2}w^{q_*}\dy\\
&+\frac{\lambda r^{\delta_*}\tau^p}{p}\int_{\R_+^N} |y|^{-s_1}w^p\dy.
\end{align*}
Choose $L>0$ so large that the sum of the first two terms is negative at $\tau=L$ and its maximum on $[0,L]$ equals the maximum on $[0,\infty)$. Then $J_{\lambda,q_*}(Lw_r)<0$ for all sufficiently small $r$. For each such fixed $r$, $J_{\lambda,q}(\tau w_r)\to J_{\lambda,q_*}(\tau w_r)$ uniformly for $\tau\in[0,L]$, as in the proof of Theorem \ref{thm:endpoint}. Hence $Lw_r$ is a negative-energy endpoint for $q$ sufficiently close to $q_*$. Maximizing the difference of the quadratic and $q_*$-power terms gives
\begin{align*}
\limsup_{q\to q_*^-}c_\lambda(q)
&\leq\max_{0\leq\tau\leq L}J_{\lambda,q_*}(\tau w_r)\\
&\leq\left(\frac12-\frac1{q_*}\right)
\left(\frac{\displaystyle\int_{\R_+^N}|\nabla w|^2\dy}
{\displaystyle(\int_{\R_+^N}|y|^{-s_2}w^{q_*}\dy)^{2/q_*}}\right)^{q_* /(q_*-2)}
+\frac{\lambda r^{\delta_*}L^p}{p}\int_{\R_+^N} |y|^{-s_1}w^p\dy.
\end{align*}
Letting $r\downarrow0$ and then taking the infimum over $w$ (nonnegative smooth functions give the same constant by approximation of the absolute value), we obtain
\begin{equation}\label{eq:upperlevel}
\limsup_{q\to q_*^-}c_\lambda(q)\leq c_+.
\end{equation}

By Theorem \ref{thm:estimate}, the energies of an arbitrary sequence of positive solutions are uniformly bounded below when $q_n\geq q_1$. If the lower limit is finite, choose a subsequence on which the energies converge to that lower limit. The same theorem bounds the subsequence in $\Hone$, and any weak limit must be zero by Lemma \ref{lem:criticalweak} and the nonexistence assumption. It follows from \eqref{eq:zerothreshold} that \eqref{eq:lowerlevel} holds. If the lower limit is $+\infty$, the inequality holds as well.

For solutions at $c_\lambda(q_n)$, the opposite upper bound follows from \eqref{eq:upperlevel}. For least-energy positive solutions it follows by Corollary \ref{cor:ground}, since their energies do not exceed $c_\lambda(q_n)$. Thus \eqref{eq:criticallevel} holds for both choices. In either case, the sequence is bounded in $\Hone$, every weakly convergent subsequence has limit zero, and the uniform exterior estimates give $u_n\to0$ in $C^1$ locally on $\overline\Omega\setminus\{0\}$. Consequently $u_n\rightharpoonup0$ in $\Hone$ and the compact weighted embedding gives $\int_\Omega u_n^p\dmu{s_1}\to0$. The energy identity in the proof of Lemma \ref{lem:criticalweak}, together with \eqref{eq:criticallevel}, now gives \eqref{eq:criticalmass}.

By \eqref{eq:criticalmass} and the Nehari identity, both measures in \eqref{eq:measure} have the required limiting total mass. Their masses outside each $B_r$ tend to zero by the convergence away from the origin. For any continuous function on $\overline\Omega$, subtract its value at $0$, split the integral into $B_r$ and its complement, and let first $n\to\infty$ and then $r\downarrow0$. This gives \eqref{eq:measure}.
\end{proof}

\section*{Acknowledgments}

This work was supported by the National Natural Science Foundation of China (No. 12601208), the Postdoctoral Fellowship Program of CPSF (No. GZC20261664), and the Innovation Research 2035 Pilot Plan of Southwest University (No. SWU-XDPY22015).

\section*{Use of generative AI tools}

During manuscript preparation, the authors used ChatGPT (OpenAI) as an auxiliary tool for literature searches, the development and checking of mathematical arguments, and editing of the exposition. The authors take full responsibility for the final manuscript.

\bibliographystyle{amsplain}
\bibliography{references}
\end{document}